\documentclass[a4paper]{article}
\usepackage[utf8]{inputenc}

\usepackage[affil-it]{authblk}
\usepackage{amsfonts, amsmath, amssymb, amsthm,fixltx2e}
\usepackage[square,sort&compress,comma,numbers]{natbib}

\usepackage{tikz}
\usetikzlibrary{shapes}
\usetikzlibrary{decorations.pathreplacing}

\usepackage{units}

\newtheorem*{conjecture*}{Conjecture}
\newtheorem{theorem}{Theorem}[section]
\newtheorem*{theorem*}{Theorem}
\newtheorem{lemma}[theorem]{Lemma}
\newtheorem{corollary}[theorem]{Corollary}
\newtheorem{proposition}[theorem]{Proposi{t}{i}on}

\theoremstyle{definition}

\theoremstyle{remark}

\newtheorem*{remark*}{Remark}

\numberwithin{equation}{section}

\newcommand{\abs}[1]{\ensuremath{\left\lvert\, #1 \,\right\rvert}}

\renewcommand{\a}{\ensuremath{\mathbf{a}}}

\newcommand{\BA}{\ensuremath{\mathrm{BA}}}
\newcommand{\C}{\ensuremath{\mathcal{C}}}
\newcommand{\comment}[1]{}
\newcommand{\CAW}{CAW}
\newcommand{\eps}{\ensuremath{\varepsilon}} 

\newcommand{\floor}[1]{\ensuremath{\left\lfloor\, #1 \,\right\rfloor}}

\renewcommand{\H}{\ensuremath{\mathcal{H}}}

\newcommand{\N}{\ensuremath{\mathbb{N}}}

\renewcommand{\P}{\ensuremath{\mathcal{P}}}
\newcommand{\betaP}{\ensuremath{\beta_{\P}}}
\newcommand{\Q}{\ensuremath{\mathbb{Q}}}
\newcommand{\R}{\ensuremath{\mathbb{R}}} 
\newcommand{\Radius}[1]{\ensuremath{\operatorname{radius} \left(\,#1\,\right)}}

\newcommand{\SP}{\ensuremath{\Sigma_{\P}}}
\newcommand{\thetaP}[1]{\ensuremath{\theta_{\P} \left(\,#1\,\right)}}
\newcommand{\x}{\ensuremath{\mathbf{x}}}
\newcommand{\XP}{\ensuremath{X_{\P}}}
\newcommand{\y}{\ensuremath{\mathbf{y}}}
\newcommand{\Z}{\ensuremath{\mathbb{Z}}}
\newcommand{\z}{\ensuremath{\mathbf{z}}}

\makeatletter
\newcommand{\subjclass}[2][1991]{%
  \let\@oldtitle\@title%
  \gdef\@title{\@oldtitle\footnotetext{#1 \emph{Mathemat{i}cs Subject Classif{i}cat{i}on.} #2.}}%
}
\newcommand{\keywords}[1]{%
  \let\@@oldtitle\@title%
  \gdef\@title{\@@oldtitle\footnotetext{\emph{Key words and phrases.} #1.}}%
}
\makeatother

\begin{document}

  \title{Cylinder absolute games on solenoids}
  \author{L.~Singhal}
  \affil{Beijing Internat{i}onal Center for Mathemat{i}cal Research,\\ Peking University, Beijing, P.~R.~China 100\,871\\
   \textit{E-mail:} {lsd\symbol{64}bicmr.pku.edu.cn}}
  
  \subjclass[2010]{Primary 11J61; Secondary 28A80, 37C45}
  \keywords{Dynamical systems, Hausdorf{f} dimension, Schmidt's game}

  \maketitle
  
  \begin{abstract}
     Let $A$ be any af{f}ine surject{i}ve endomorphism of a solenoid \SP\ over the circle $S^1$ which is not an inf{i}nite-order translat{i}on of \SP. We prove the existence of a cylinder absolute winning (\CAW) subset $F \subset \SP$ with the property that for any $x \in F$, the orbit closure $\overline{\{ A^{\ell} x \mid \ell \in \N \}}$ does not contain any periodic orbits. The class of inf{i}nite solenoids considered in this paper provides, to our knowledge, some of the f{i}rst examples of non-Federer spaces where absolute games can be played and won. Dimension maximality and incompressibility of \CAW\ sets is also discussed for a number of possibilit{i}es in addit{i}on to their winning nature for the games known from before.
  \end{abstract}

  \section{Introduction}\label{S:intro}
    Let \P\ be an (f{i}nite or inf{i}nite) ordered set of prime numbers with $p_1 < p_2 < \cdots$ and $X_{\P, n}$ be the restricted product space
    \begin{equation}\label{E:spaceX}
      \R^n \times \prod{}^{\prime}_{p \in \P}\, \Q_{p}^n,
    \end{equation}
    where $\prod'$ denotes that for each element $x \in X_{\P, n}$, the entries $x_p \in \Z^n_p$ for all but f{i}nitely many $p$'s. A \emph{\P-solenoid} of topological dimension $n$ is the quo{t}{i}ent space $\SP (n) := \nicefrac{X_{\P, n}}{\Delta ( R^n )}$, where $R$ is the ring $\Z [ \{ p^{-1} \mid p \in \P \} ]$ whose $n$-fold product is embedded diagonally in $X_{\P, n}$ as a uniform la{t}{t}{i}ce. We call the quo{t}{i}ent map $X_{\P, n} \rightarrow \SP (n)$ to be $\Pi$. Solenoids are compact, connected metrizable abelian groups. They have somet{i}mes been called ``fractal versions of tori''~\citep{Sem12c}. When \P\ is a f{i}nite set of cardinality $l - 1$, the Hausdorf{f} dimension of $X_{\P, n}$ (and therefore of $\SP (n)$ too) under the natural metric given by~\eqref{E:solmet} is $nl$. This also implies that the dimension is inf{i}nite when \P\ is so, as the increasing sequence of f{i}nite products associated with the f{i}nite truncat{i}ons of \P\ are isometrically embedded inside $X_{\P, n}$.\\[-0.1cm]

    The set of endomorphisms of $\SP := \SP (1)$ is precisely the ring $R$ whose elements act mul{t}{i}plica{t}{i}vely componentwise. An \emph{af{f}ine transformat{i}on} $A : \SP \rightarrow \SP$ is meant to denote the map
    \begin{equation}\label{E:affine}
      \mathbf{x} \mapsto ( m / n ) \mathbf{x} + \mathbf{a}
    \end{equation}
    where $m / n \in R \setminus \{ 0 \}$ and $\mathbf{a} \in \SP$. It is well known that when $A = m /n$ is a surject{i}ve endomorphism of the solenoid, it acts ergodically on \SP\ if{f} $m/n \notin \{ 0, \pm 1 \}$~\citep[Proposit{i}on~1.4]{Wil76}. We also learn from \citeauthor{Ber85}~\citep[Theorem~3.2]{Ber85} that for every compact group $G$, any semigroup of its af{f}ine transformat{i}ons lying above an ergodic semigroup of surject{i}ve endomorphisms is ergodic as well. This gives us a suf{f}icient condit{i}on for the transformat{i}on $A \mathbf{x} = ( m / n )\x + \mathbf{a}$ to be ergodic, namely that $m/n \neq \pm 1$. Ergodicity of the act{i}on guarantees that almost all orbits of $A$ are dense in \SP. However, just like \citeauthor{Dan88}~\citep{Dan88}, this work is concerned with understanding the complementary set. Given an af{f}ine transformat{i}on $A$, we will like to know the set of points of \SP\ whose $A$-orbits remain away from periodic $B$-orbits for all $B \in R \setminus \{ \pm 1 \}$.\\[-0.1cm]
    
    When \P\ is the set consist{i}ng of all the primes in \N, the space \SP\ is called the \emph{full solenoid} over $S^1$ with the f{i}eld \Q\ being its ring of endomorphisms. Let $B$ be a non-zero rat{i}onal number. The growth of the number of $B$-periodic orbits as a funct{i}on of the period is determined by the entropy of the act{i}on on \SP. Lat{t}er has been computed by \citeauthor{Yuz}~\citep{Yuz} f{i}rst and recovered by \citeauthor{LW88} in~\citep{LW88} who explained it to be the sum of the Euclidean and the $p$-adic contribut{i}ons. In fact, they achieve it for all automorphisms of solenoids over higher-dimensional tori as well. We remark that each such epimorphism li{f}{t}s uniquely to a homomorphism from $\XP := X_{\P, 1}$ to itself, which we shall con{t}{i}nue to denote by the same ra{t}{i}onal number. For an af{f}ine transformat{i}on, we however have a choice involved in terms of a representat{i}ve for the translat{i}on part $\mathbf{a}$.\\[-0.1cm]
  
    Let $y \in \SP$ be arbitrary and $A$ be a (surject{i}ve) af{f}ine transformat{i}on of \SP\ as in~\eqref{E:affine} with either $m/n \neq 1$ or $\mathbf{a} = \mathbf{0}$. We intend to show that the set of points $\mathbf{x} \in \XP$ whose forward orbit under the map $\mathbf{x} \mapsto A\mathbf{x}$ maintains some posi{t}{i}ve distance from the $1$-uniformly discrete subset $\Pi^{-1} ( \{ y \} ) \subset \XP$ is \emph{cylinder absolute winning} (\CAW) in a similar sense as \citeauthor*{FSU15}~\citep{FSU15}. Once we have this, we can take intersect{i}on of countably many of these sets to conclude about $A$-orbits which avoid neighbourhoods of all periodic orbits of surject{i}ve endomorphisms. This strategy is in same taste and builds upon the work of \citeauthor{Dan88}~\citep{Dan88} on orbits of semisimple toral automorphisms.\\[-0.1cm]
    
    Our setup has two players in which one of them (Alice) will be blocking open cylinder subsets of \XP\ at every stage of a two-player game. To elaborate, one such cylinder is given by
    \begin{equation}\label{E:cyl}
      C ( \mathbf{x}, \eps, i ) := \begin{cases}
        \R \times \prod_{j < i} \Q_{p_j} \times B ( x_i, p_i \eps ) \times \prod{}^{\prime}_{j > i} \Q_{p_{j}} &\textrm{ if } i > 0 \textrm{ and}\\
        B ( x_0, \eps ) \times \prod^{\prime}_{j > 0} \Q_{p_j} &\textrm{ otherwise,}
      \end{cases}
    \end{equation}
    where $\mathbf{x} = ( x_0, \ldots, x_i, \ldots ) \in \XP$. For us, $B ( x_i, r )$ will be always be the set of points in $\Q_{p_i}$ whose distance from $x_i$ is strictly less than $r$ while $\overline{B} ( x_i, r )$ will also include those whose distance from $x_i$ is exactly $r$. We explain this game in \S\,\ref{S:games} a{f}{t}er a brief tour of some of its older and related versions. Our aim is to prove the following statement in this paper:
    {\begin{theorem}\label{Th:main}
      Let \SP\ be a solenoid over the circle $S^1$ and
      \[
        \{ A_j : \x \mapsto ( m_j / n_j ) \x + \a_j \mid j \in \N\}
      \]
      be any subset of af{f}ine surject{i}ve endomorphisms of \SP\ such that
      \begin{enumerate}
        \item none of the $A_j$'s is a non-trivial translat{i}on of \SP, and
        \item the collect{i}on of rat{i}onal numbers $\{ m_j / n_j \}_{j \in \N}$ lying below the family $\{ A_j \}$ belong to some f{i}nite ring extension of \Z.
      \end{enumerate}
       Then, there exists a cylinder absolute winning subset $F \subset \SP$ such that for any $x \in F$ and $j \in \N$, the orbit closure $\overline{ \{ A^{k}_{j} x \mid k \in \N \} }$ contains no periodic $B$-orbit for all $B \in R \setminus \{ \pm 1 \}$.
    \end{theorem}}
    \noindent {This is done in \S\,\ref{S:morph}. In the last sect{i}on, we illustrate how this informat{i}on can be used to infer something about the $f$-dimensional Hausdorf{f} measure of $F$. We also discuss the strong game winning and incompressible nature of \CAW\ subsets when \P\ is f{i}nite.}\\[-0.1cm]
    
    \subsection{\protect{Comparison with the work of \citeauthor{Wei13}}}\label{SS:Weil}
      It is plausible that some of our results given here may also be obtained from the more general framework discussed by~\citeauthor{Wei13}~\citep{Wei13}. We take some t{i}me to expound the import of his main result.\\[-0.1cm]
    
      Let $( \overline{X}, d )$ be a proper metric space (i.~e.\ all closed balls are compact) and $X$ a closed subset of $\overline{X}$. Consider a family of subsets $\{ R_{\lambda} \subset \overline{X} \mid \lambda \in \Lambda \}$ which are called \comment{resonant sets}\emph{resonant sets} and a family of contract{i}ons $\{ \psi_{\lambda} \mid\ ] 0, 1 ] \rightarrow 2^{\overline{X}} \mid \lambda \in \Lambda \}$, indexed by some (same) countable set $\Lambda$. It is required that $R_{\lambda} \subset \psi_{\lambda} ( t_1 ) \subset \psi (t_2)$ for all $0 < t_1 < t_2$ and $\lambda$. This datae is writ{t}en in a concise form as $\mathcal{F} = ( \Lambda, R_{\lambda}, \psi_{\lambda} )$. The set of \comment{approximable! badly}\emph{badly approximable points in $S$} with respect to the family $\mathcal{F}$ is def{i}ned as
      \begin{equation}\label{E:BAF}
        \BA_X (\mathcal{F}) := \big\{ x \in X \mid \exists c = c(x) > 0 \textrm{ such that } x \notin \bigcup_{\lambda \in \Lambda} \psi_{\lambda} (c) \big\}.
      \end{equation}
      Next, each $R_{\lambda}$ is assigned a \comment{height funct{i}on}\emph{height} $h_{\lambda}$ with $\inf_{\lambda} h_{\lambda} > 0$. The standard contract{i}on $\psi_{\lambda}$ is then determined as $\psi_{\lambda} (c) := \mathcal{N}_{c / h_{\lambda}} (R_{\lambda})$, where $\mathcal{N}_{\eps} (S)$ denotes the set of all points of $\overline{X}$ in the $\eps$-vicinity of elements of $S$. We further assume that our resonant sets $R_{\lambda}$ are \comment{resonant sets! nested}\emph{nested} with respect to the height funct{i}on $h$, i.~e., $R_{\lambda} \subseteq R_{\beta}$ for every $\lambda, \beta \in \Lambda$ such that $h_{\lambda} \leq h_{\beta}$ and that the values taken by $h$ form a discrete subset of\ $]\,0,\,\infty\,[$. For any collect{i}on $\mathcal{S}$ of subsets of $\overline{X}$, the set $X$ is said to be \comment{$b_*$-di{f}fuse}\emph{$b_*$-di{f}fuse with respect to $\mathcal{S}$} for some $0 < b_* < 1$ if there exists some $r_0 > 0$ such that for all balls $\overline{B} ( x, r ),\ x \in X,\ 0 < r < r_0$ and $S \in \mathcal{S}$, there exists a sub-ball
      \begin{equation}\label{E:diffuse}
        \overline{B} ( y, b_* r ) \subset \overline{B} ( x, r ) \setminus \mathcal{N}_{b_* r} (S) \textrm{ with } y \in X.
      \end{equation}
      The family $\mathcal{F}$ is \comment{locally contained}\emph{locally contained in $\mathcal{S}$} if for any $\overline{B} ( x, r ),\ x \in X$ with $r < r_0$ and $\lambda \in \Lambda$ with $h_{\lambda} \leq 1/r$, there is an $S \in \mathcal{S}$ such that $B \cap R_{\lambda} \subset S$.
      Concrete realizat{i}ons of this abstract formalism include the case when $\overline{X}$ is the Euclidean space $\R^n$ and $\mathcal{S}$ consists of all af{f}ine hyperplanes in $\overline{X}$. Some examples of \emph{hyperplane dif{f}use} sets are supports of absolutely decaying measures such as the Cantor sets and the Sierpi\'nski triangle.
     \begin{theorem}[\protect{\citeauthor{Wei13},~\citeyear{Wei13}}]
        Let $X \subset \overline{X}$ be closed and $b_*$-dif{f}use with respect to a collect{i}on $\mathcal{S}$ of subsets of $\overline{X}$. Also, $\mathcal{F}$ is a family with nested resonant sets $R_{\lambda}$ and discrete heights, locally contained in $\mathcal{S}$. The set $\BA_X (\mathcal{F})$ def{i}ned in~\eqref{E:BAF} is Schmidt winning subset of $X$.
      \end{theorem}    
      \noindent Many of the terms used above will be explained in \S~\ref{S:games}. Actually, it is shown in~\citep{Wei13} that $\BA_X (\mathcal{F})$ is \emph{absolute winning with respect to $\mathcal{S}$}. This covers many important examples like the set of \comment{approximable! badly}$( s_1, \ldots, s_n )$-badly approximable vectors in $\R^n,\ \C^2$ and $\Z_p^2$, the set of sequences in the Bernoulli-shift which avoid all periodic sequences and the set of orbits of toral endomorphisms which stay away from periodic orbits. Curiously for us, we do not f{i}nd any discussion on solenoidal endomorphisms in his work.\\[-0.1cm]
    
      In order to be able to use his results, we would have to show that (a) the space \XP\ is $b_*$-dif{f}use with respect to the collect{i}on $\mathcal{C}$ of open cylinders in \XP\ for an appropriate value of $b_*$, and (b) the family of pre-images $A^{-j} B ( \y + \z, t )$ is locally contained in $\mathcal{C}$. We have endeavoured to provide a more direct proof here. In this sense, our work may be considered as an addit{i}on to the list of examples given in~\citep{Wei13}. We must also point out that the assumpt{i}on of having a Federer measure with $X$ as its support in~\citep[Proposit{i}on~2.1]{Wei13} need not be true when $X = \XP,\ \P$ is inf{i}nite and $\mu$ is the Haar measure on \XP\ (cf.~\eqref{E:nonFed}).

  \section{Metric and measure structure on solenoids}\label{S:sol}
    Before we discuss the game, it is imperat{i}ve that we say a few words about how balls in our space $X_{\P, n}$ ``look like.'' A metric on $X_{\P, n}$ is given by
    \begin{equation}\label{E:solmet}
      d ( \mathbf{x}, \mathbf{z} ) := \max \left\{ \abs{x_0 - z_0},\ \sup_{p \in \P} \big\{ p^{-1} \abs{x_p - z_p}_{p} \big\} \right\}
    \end{equation}
    where \abs{\cdot}\ is the usual Euclidean metric on $\R^n$ and $\abs{\cdot}_{p}$ refers to the $p$-adic ultrametric on $\Q_{p}^n$ such that the diameter of $( p^{-1} \Z_{p} )^n$ is $p$. Clearly, distance between any two dis{t}{i}nct points of $\Delta (R^n)$ is at least $1$ or in other words, the injec{t}{i}vity radius for the quo{t}{i}ent map $\Pi : X_{\P, n} \rightarrow \SP (n)$ equals $1$. More generally, as per~\citep{BFK11}, a subset $Z$ of any metric space $X$ is said to be \emph{$\delta$-uniformly discrete} if the distance between any two distinct points of $Z$ is at least $\delta$. In this terminology, the set $\Delta ( R^n )$ described above is $1$-uniformly discrete in $X_{\P, n}$.\\
    
    The definit{i}on of the metric in~\eqref{E:solmet} ensures that balls $B ( \x, r ) \subset \XP$ are direct products of their coordinatewise project{i}ons when (i) \P\ is a f{i}nite set, or (ii) when $r < 1$. Moreover for $r < 1,\ 1 \leq p_i r < p_i$ for all $i$ large enough. As the coordinates $x_i$ of \x\ belong to $\Z_{p_i}$ for $i \gg 1$, the project{i}ons $B ( x_i, p_i r ) = \Z_{p_i}$ for all but f{i}nitely many $i$'s whereby we get that any such `open' ball $B ( \x, r )$ is actually open in \XP\ while the `closed' ball $\overline{B} ( \x, r )$ is a compact neighbourhood of \x. It is easy to see that the diameter of $B ( \x, r )$ is $2r$ but in most (all but one) of the places $p \in \P$, the $p$-adic diameters of the project{i}ons will be strictly less than $pr$ as distances in the non-archimedean f{i}elds are a discrete set. For a given $r > 0$, there might not be any integer $j$ such that $p_i^j = p_i r$. This seemingly minor issue is very crucial in our next set of calculat{i}ons. If $j_i \in \Z$ is such that $p_i^{j_i} \leq r < p_i^{j_i + 1}$, then we call $\lfloor r \rfloor_i := p_i^{j_i}$. It should be noted that $\floor{p_i^m r}_i = p_i^m \floor{r}_i$ for all $m \in \Z$ and more generally, $\floor{t r}_i \leq \floor{p_i \floor{t}_i r}_i = p_i \floor{t}_i \floor{r}_i\ \forall r, t > 0$. The lemma below may be of independent interest to the reader:
    \begin{lemma}\label{L:prod}
      Let $x > 1$ and $P_{\P} (x)$ denote the product
      \begin{equation}
        \prod_{i \in \N} \floor{\frac{p_i}{x}}_{i} = \prod_{\substack{i \in \N,\\ p_i < x}} \floor{\frac{p_i}{x}}_{i}.
      \end{equation}
      Then,
      \begin{equation}
        \ln P_{\P} (x) \leq - \thetaP{x} + O \big( ( \ln x )^2 \big), \textrm{ where } \thetaP{x} := \sum_{\substack{p \leq x\\ p \in \P}} \ln p.
      \end{equation}
    \end{lemma}
    \begin{proof}
      Let $p_i < x$ and $m_i \in \N$ be the unique integer for which
      \begin{equation}
        \dfrac{1}{p_i^{m_i}} \leq \dfrac{p_i}{x} < \dfrac{1}{p_i^{m_i - 1}}\quad \Leftrightarrow\quad p_i^{m_i} < x \leq p_i^{m_i + 1}.
      \end{equation}
      We know that $m_i = k$ if and only if $x^{1 / ( k + 1 )} \leq p_i < x^{1/k}$. This observat{i}on leads to decomposing $P_{\P} (x)$ as a double product
      \begin{equation}
        P_{\P} (x) = \prod_{k = 1}^{\ell} \prod_{\substack{p \in \P,\\ x^{1 / ( k + 1 )} \leq p < x^{1/k}}} p^{-k}
      \end{equation}
      where $\ell \in \N$ is such that $x^{1 / ( \ell + 1 )} \leq 2 < x^{1 / \ell}$ (i.~e., $\log x / \log 2 - 1 \leq \ell < \log x / \log 2$). Taking negat{i}ve logarithms on both sides, we have
      \begin{align}
        - \ln P_{\P} (x) &= \sum_{k = 1}^{\ell} k \sum_{\substack{p \in \P,\\ x^{1 / ( k + 1 )} \leq p < x^{1/k}}} \ln p\notag\\
                         &\geq \sum_{k = 1}^{\ell} k \left( \thetaP{x^{1/k}} - \thetaP{x^{1/ ( k + 1 )}} - \frac{1}{k} \ln x \right)\\
                         &= \sum_{k = 1}^{\ell} \thetaP{x^{1/k}} - \ell \left( \thetaP{x^{1/ ( \ell + 1 )}} + \ln x \right).\notag
      \end{align}
      As $\thetaP{x^{1 / ( \ell + 1 )}} \leq \thetaP{2} \leq \ln 2$ and $\ell \ll \ln x$, we are done.
    \end{proof}
    \noindent When \P\ is the full subset of primes, $\thetaP{x} \sim x$ as $x \rightarrow \infty$. If \P\ consists of all primes in an arithmet{i}c progression with common dif{f}erence $k$ and $( a, k ) = 1$ where $a$ is one of the terms in the progression, then \thetaP{x} goes roughly as $x / \varphi (k)$ as $x \rightarrow \infty$. Here, $\varphi$ stands for the Euler's tot{i}ent funct{i}on. For more on this, the reader is redirected to \citeauthor{MV07}'s book~\citep{MV07}. We conclude that for all subsets \P\ of primes which come from arithmet{i}c progressions (and in part{i}cular the full subset), there exist some $0 < c_1 = c_1 (k) < 1$ such that
    \begin{equation}\label{E:Pp1byr}
      P_{\P} \left( \frac{1}{r} \right) \leq r^{c_1 \frac{1/r}{\ln (1/r)}} \textrm{ for all } 0 < r \ll 1.
    \end{equation}
    This short exercise serves dual purpose. On one hand, it is a roundabout way of establishing the inf{i}niteness of the Hausdorf{f} dimension of \XP\ when \P\ is as above using the \emph{mass distribut{i}on principle}. A f{i}nite non-zero measure $\mu$ whose support is a bounded subset of a metric space $M$ is called a \emph{mass distribut{i}on} on $M$. We require our \emph{dimension funct{i}ons} (also known as \emph{gauge funct{i}ons}) $f : \R_{\geq 0} \rightarrow \R_{\geq 0}$ to be increasing in some neighbourhood $[ 0, a_f )$, cont{i}nuous on $( 0, a_f )$, right cont{i}nuous at $0$ and $f (r) = 0$ if{f} $r$ equals zero~\citep[pg.~33]{Fal90}.
    \begin{proposition}[cf.~\citep{Fal90}]\label{P:massdis}
      Let $\mu$ be a mass distribut{i}on on a second countable metric space $M$ such that for some dimension funct{i}on $f$ and $\delta_0 > 0$,
      \[
        \mu (U) \leq c_2 f ( \abs{U} )
      \]
      for some f{i}xed $c_2 > 0$ and all subsets $U$ with $\abs{U} \leq \delta_0$. Then, $\mathcal{H}^f (M) \geq \mu (M) / c_2$. If $f$ is the power rule $r \mapsto r^s$, we can say that $\dim M \geq s$.
    \end{proposition}
    Now, let $\mu$ be the restrict{i}on to $[ 0, 1 ] \times \prod_{p \in \P} \Z_p$ of the Haar measure $\nu$ on \XP\ which is the product of the Haar measures on \R\ and on each $\Q_p,\ p \in \P$. We normalize it so that $\mu$ is a probability measure on \XP. Any ball $B ( \x, r )$ with radius $0 < r \ll 1$ will then have
    \begin{equation}\label{E:nonFed}
      \mu \big( B ( \x, r ) \big) \leq 2r \times \prod_{i \in \N} \lfloor p_i r \rfloor_i = 2r \cdot P_{\P} ( 1/r ) \leq 2 r^{c_1 \frac{1/r}{\ln ( 1/r )} + 1}
    \end{equation}
    by~\eqref{E:Pp1byr} when \P\ is the set of all primes in any inf{i}nite arithmet{i}c progression with $( a, k ) = 1$. Af{t}er some more work, the above proposit{i}on will give us that the Hausdorf{f} dimension of the space \XP\ is inf{i}nite in such cases.\\[-0.1cm]

    On the other hand, Lemma~\ref{L:prod} will help us again in \S\,\ref{S:hausdorff} when we examine the dimension-theoret{i}c largeness of various winning subsets of cylinder absolute games. In our version, Alice shall be dealing with the family \C\ of closed subsets exactly one of whose co-ordinates $x_i,\ i \in \{ 0, \ldots, l - 1 \}$ is a fixed constant. The resulting $\eps$-neighbourhood $C ( \x, \eps, i )$ of such a set $P \in \C$ as also de{f}{i}ned in~\eqref{E:cyl} will be called an (open) \emph{cylinder}. Given a cylinder $C = C ( \x, \eps, i )$, we say that the \emph{radius of $C$} is $\eps$ if $i = 0$ and the minimum such $\eps'$ for which $C ( \x, \eps', i ) = C ( \x, \eps, i )$ otherwise. The index $i$ is called the \emph{constraining coordinate of $C$}. We emphasize that both Alice and Bob are fully aware of the radii of the balls chosen by the lat{t}er at any stage of our game by reading the real coordinate.
    
  \section{In{f}{i}nite games on complete metric spaces}\label{S:games}
    Let $M$ be a complete metric space and $F$ be a {f}{i}xed subset of $M$. In the original game introduced by	\citeauthor{Sch66}~\citep{Sch66}, Alice and Bob are two players who each take turns to pick closed balls in $M$ in the following manner: We have $\alpha, \beta \in ( 0, 1 )$ to be two real numbers such that $1 - 2 \alpha + \alpha \beta > 0$. The game begins with Bob choosing any closed ball $B_0 = \overline{B} ( b_0, r ) \subseteq M$ subsequent to which Alice has to make a choice of some $A_1 = \overline{B} ( a_1, \alpha r )$ such that $A_1 \subset B_0$. After this, Bob picks $B_1 = \overline{B} ( b_1, \beta \alpha r ) \subset A_1$ and the game goes on till in{f}{i}nity. We thus get a decreasing sequence of closed, non-empty subsets of a complete metric space
    \begin{equation}\label{E:seq}
      M \supset B_0 \supset A_1 \supset B_1 \supset A_2 \supset \ldots
    \end{equation}
    Alice is declared the winner if $\cap_j B_j = \cap_j A_j = \{ a_{\infty} \} \subset F$. A set $F$ is called \emph{$( \alpha, \beta )$\,-\,winning} if Alice has a strategy to win the above game regardless of Bob's moves. Further, it is \emph{$\alpha$\,-\,winning} if it is $( \alpha, \beta )$\,-\,winning for all $\beta \in ( 0, 1 )$ and \emph{Schmidt winning} if it is $\alpha$-winning for some $\alpha \in ( 0, 1 )$.\\[-0.1cm]
  
    For various applica{t}{i}ons of pract{i}cal interest, one {f}{i}nds out that Alice need not bother herself too much about choosing the balls $A_j$'s as long as she is able to block out neighbourhoods of certain undesirable points. This is true for example when $M$ is the real line and $F$ is the set BA of badly approximable numbers as discussed in~\cite{Sch66} where Alice needs to be far from ra{t}{i}onal numbers with small denominators. Moreover if she is careful enough about her strategy, she has to worry about very few of such ra{t}{i}onals -- at t{i}mes just one of them and hence, she need only shift the game outside of a ball B centered at some $p/q \in \Q \cap B_j$ at the $j$-th stage. This was formalized by \citeauthor{McM10}~\citep{McM10} who called the new variant to be \emph{absolute winning games}. Let $\beta \in ( 0, 1/3 )$ and now Alice chooses open balls $A_j$'s with radius $(A_j) \leq \beta \cdot \mathrm{radius}\,(B_{j - 1})$. Then, Bob is supposed to pick some closed ball $B_j \subset B_{j - 1} \setminus A_j$ with $\mathrm{radius}\,( B_j ) \geq \beta \cdot \mathrm{radius}\,( B_{j - 1} )$. The sequence of sets we now have is
    \begin{equation}\label{E:abs}
      B_0 \supset B_0 \setminus A_1 \supset B_1 \supset B_1 \setminus A_2 \supset \ldots
    \end{equation}
    Note that the countable intersec{t}{i}on $\cap_j B_j$ might be bigger than a singleton set now and we have to set the winning condi{t}{i}on to be $\cap_j B_j \cap F \neq \phi$. The set $F$ is said to be \emph{$\beta$-absolute winning} if Alice has a strategy to win in this situa{t}{i}on and it is called \emph{absolute winning} if it is $\beta$-absolute winning for all $\beta \in ( 0, 1/3 )$.\\[-0.1cm]
    
    In the same paper~\citep{McM10}, \citeauthor{McM10} also gave the concept of a \emph{strong winning set.} We again have two parameters $\alpha, \beta \in\,]\,0,\,1\,[$ but now Alice is allowed more opt{i}ons in the form of balls $A_{i + 1} \subset B_i$ such that $\abs{A_{i + 1}} \geq \alpha \abs{B_i}$ for all $i \in \N$ while for Bob, $\abs{B_i} \geq \beta \abs{A_i}$ for $i > 1$. It cont{i}nues to be mandatory that $A_i \subset B_{i - 1}$ and $B_i \subset A_i$ for all $i$. A subset $F$ for which Alice has a winning strategy in this game is called an \emph{$( \alpha, \beta )$-strong winning set.} The subset $F$ is said to be \emph{$\alpha$-strong winning} if it is $( \alpha, \beta )$-strong winning for all $\beta \in\,]\,0,\,1\,[$ and \emph{strong winning} if it is $\alpha$-strong winning for some $\alpha > 0$. For Euclidean spaces, a strong winning subset is Schmidt winning too and retains its strong winning property under quasisymmetric mappings~\citep[Theorem~1.2]{McM10}.\\[-0.1cm]
    
    The absolute game has an obvious drawback that if $F$ is the set of badly approximable vectors in $\R^n$ for any $n > 1$, then Bob can force the game to be always centered on the hyperplane $\R^{n - 1} \times \{ 0 \}$ and Alice is not able to win trivially. Therefore, it was proposed in~\citep{BFK+12} that she be allowed to block out a neighbourhood of some $k$\,-\,dimensional a{f}{f}{i}ne subspace of $\R^n$ at each stage of the game. Taking this into considera{t}{i}on, they gave a family of games played on the Euclidean space $\R^n$ called \emph{$k$-dimensional $\beta$-absolute games} ($0 < \beta < 1/3,\ 0 \leq k < n$) where Bob having chosen $B_0 = B ( \mathbf{b}_0, r_0 ) \subset \R^n$, Alice picks some a{f}{f}{i}ne subspace $V_1$ of dimension $k$ and for some $0 < \varepsilon_1 \leq \beta r_1$ removes the $\varepsilon_1$-neighbourhood of $V_1$, namely $A_1 = V_1^{(\varepsilon_1)}$ from $B_0$. This is followed by Bob picking a closed ball $B_1 \subseteq B_0 \setminus A_1$ with $\mathrm{radius}\,(B_1) \geq \beta r_0$ and the game proceeds in a similar fashion. In general, the parameter $\varepsilon_j$ is allowed to depend on $j$ subject only to $0 < \varepsilon_j \leq \beta\cdot\mathrm{radius} (B_j)$. Alice wins if $\cap_j B_j \cap F \neq \phi$. As before, $F \subseteq \R^n$ is \emph{$k$-dimensional $\beta$-absolute winning} if Alice can win the $k$-dimensional $\beta$-absolute game over $F$ irrespec{t}{i}ve of Bob's strategy. It is called \emph{$k$-dimensional absolute winning} if it is $k$-dimensional $\beta$-absolute winning for all $\beta \in ( 0, 1/3 )$. It is clear from the de{f}{i}ni{t}{i}ons that for $0 \leq k_1 < k_2 < n$, if a set $F \subseteq \R^n$ is $k_1$-dimensional $\beta$-absolute winning, then it is $k_2$-dimensional $\beta$-absolute winning too. Also, $0$-dimensional $\beta$-absolute winning is the same as $\beta$-absolute winning.\\[-0.1cm]
    
  All of this culminated in the axioma{t}{i}za{t}{i}on by \citeauthor*{FSU15}~\citep{FSU15} where $M$ is a complete metric space, \H\ is a non-empty collec{t}{i}on of closed subsets of $M$ and $F \subseteq M$ is {f}{i}xed before the start of play. For $0 < \beta < 1$, the set $F$ is called \emph{$( \H, \beta )$-absolute winning} if Alice can ensure the intersec{t}{i}on $\cap_j B_j \cap F \neq \phi$ by removing neighbourhoods $A_j = H_j^{(\varepsilon_j)}$ for some $H_j \in \H$ and $0 < \varepsilon_j \leq \beta\cdot\mathrm{radius}\,(B_{j - 1})$ at every $j$-th stage of the game. We follow \citeauthor{KL15}~\citep{KL15} to declare Bob the winner by default if at any (f{i}nite) stage of the game, he is le{f}{t} with no legal choice of the ball $B_j$ to make. In the course of the game, Alice will have to make sure that such an event does not ever occur. This is keeping in mind the example of a Schmidt game illustrated in~\citep[Proposit{i}on~5.2]{KW10} where Bob is not able to win because he has no opt{i}on of $B_j$ left. The reader is caut{i}oned at this point that in \citep[Def{i}nit{i}on~C.1]{FSU15}, the authors resort to the opposite convent{i}on of Alice winning the game if it ends abruptly.\\[-0.1cm]
  
  Ever since \citep{Sch66} came out, Schmidt games have been played and won over subsets of various metric spaces. We were unable to f{i}nd any reasonable survey art{i}cle covering the developments in the area. It will also be impossible to give here a comprehensive account of all the progress that has been made by dif{f}erent people and groups. We will have to contend ourselves by pointing to only a few representat{i}ve works. \citeauthor{Dan86}~\citep{Dan86} formulated and proved results about the winning nature of the set of points in a \comment{homogeneous space}homogeneous space $G / \Gamma$ of a semisimple Lie group $G$ whose orbits under a one-parameter subgroup act{i}on are bounded. \comment{Aravinda, C.~S.}\citeauthor{Ara94}~\citep{Ara94} showed that the set of points on any non-constant $C^1$ curve $\sigma$ on the unit tangent sphere $S_p$ of any point $p$ on a complete, non-compact Riemannian manifold $M$ with constant negat{i}ve curvature and f{i}nite Riemannian volume which lead to bounded geodesic orbits is Schmidt winning.\\[-0.1cm]
  
  When $\Gamma \subset G$ is an irreducible lat{t}ice of a connected, semisimple $G$ with no compact factors, \citeauthor{KM96}~\citep{KM96} established that the subset of points in $G / \Gamma$ with bounded $H$-orbits is of full Hausdorff dimension whenever $H$ is a nonquasiunipotent one-parameter subgroup of $G$. \citeauthor{KW10}~\citep{KW10} allowed for Alice's and Bob's choices of subsets to be more flexible than just metric balls and used this to set{t}le that the set of $\mathbf{s}$-badly approximable vectors in $\R^n$ is winning for any f{i}xed $\mathbf{s} \in \R^n_+$. This was part of an ef{f}ort to understand Schmidt's conjecture on the intersect{i}on of the sets of weighted badly approximable vectors for dif{f}erent weights which was f{i}nally resolved by \citeauthor*{BPV11}~\citep{BPV11}. It was shown by \citeauthor*{EGL}~\citep{EGL} that the set of points on any $C^1$ curve which are badly approximable by rat{i}onals coming from a number f{i}eld $\mathbb{K}$ is Schmidt winning. {More generally, it is possible to def{i}ne a hyperplane absolute game on any $C^1$ manifold. For example, it was recently proved in~\cite{AGGL} that for any one parameter $\mathrm{Ad}$-semisimple subsemigroup $\{ g_t \}_{t \geq 0}$ of the product $G$ of f{i}nitely many copies of $\mathrm{SL}_2 (\R)$'s, the set of points $x$ belonging to any lat{t}ice quot{i}ent $G / \Gamma$ of $G$ and with bounded $\{ g_t \}$-orbit in $G / \Gamma$ is hyperplane absolute winning.}\\[-0.1cm]
    
  In our se{t}{t}{i}ng, $M$ shall be \XP\ (or $\Sigma_{\P}$ if you prefer), \H\ is the family \C\ of subsets described in \S\,\ref{S:sol} and an example of the target set $F$ is given below. A less contrived one will be available in the next sec{t}{i}on. A \emph{cylinder $\beta$-absolute game} begins with Bob choosing a closed ball $B_0 = \overline{B} ( \mathbf{x}_0, r_0 )$. Subsequent to this, Alice blocks an open cylinder $C_1$ whose radius has to be less than or equal to $\beta r_0$. The cylinders seem to us to be the appropriate replacement for the hyperplane neighbourhoods of~\citep{BFK+12} in metric spaces like solenoids. Recall that the exact value of the radius can be read of{f} from the real coordinate. Moreover, $\mathrm{radius}\,(B_j)$ is required to be at least $\beta\cdot\mathrm{radius}\,(B_{j - 1})$ for all $j \in \N$. The game of our interest goes as
    \begin{equation}\label{E:caw}
      B_0 \supset B_0 \setminus C_1 \supset B_1 \supset B_1 \setminus C_2 \supset \cdots
    \end{equation}
    and $F$ is said to be \emph{cylinder $\beta$-absolute winning} if Alice can devise a method to win this game, i.~e., $\bigcap_j B_j \cap F \neq \phi$. It is {\emph{cylinder absolute winning}} if there exists a $0 < \betaP = \betaP (F) \leq 1/3$ such that $F$ is cylinder $\beta$-absolute winning for all $\beta \in\ ]\,0,\,\betaP\,[$. The supremum of such $\betaP$'s is christened the \emph{\CAW\ dimension} of $F$.
    \begin{proposition}\label{P:ccwin}
      A countable intersect{i}on of \CAW\ subsets of \XP\ with winning dimension $\geq \beta_0$ each is a \CAW\ set with winning dimension $\geq \beta_0$.
    \end{proposition}
    The basic idea of the proof remains the same as in \citep[Theorem~2]{Sch66} and is being skipped here. We next give a theorem largely inspired by one of \citeauthor{Dan86}~\citep{Dan86}.
    \begin{theorem}\label{Th:cgame}
      Let $N$ be a countable indexing set and $\{ A_{( n, t )} \subseteq C_{( n, t )} \subset \XP \mid n \in N,\ t \in ( 0, 1 ) \}$ be a family of set pairs where $C_{( n, t )}$'s are restricted to be open cylinders in \XP\ with the same f{i}xed constraining coordinate $i$. If for any compact $K \subset \XP$ and $\mu \in ( 0, 1 )$, there exist $R \geq 1,\ \eps \in ( 0, 1 )$ and a sequence $( R_n )$ of posit{i}ve reals with the following propert{i}es:
      \begin{enumerate}
        \item if $n \in N$ and $t \in ( 0, \eps )$ are such that $A_{( n, t )} \cap K \neq \phi$, then $R_n \leq R$ and the radius $r ( C_{( n, t )} )$ of the cylinder $C_{(n, t )}$ is at most $t R_n$,
        
        \item if $n_1, n_2 \in N$ and $t \in ( 0, \eps )$ are such that both $A_{( n_i, t )}$ intersect $K$ non-trivially and the radius bounds of the associated cylinders are comparable, i.~e., $\mu R_{n_1} \leq R_{n_2} \leq \mu^{-1} R_{n_1}$, then either $n_1 = n_2$ or $d\,( A_{( n_1, t )},\ A_{( n_2, t )} ) \geq \eps ( R_{n_1} + R_{n_2} )$.
      \end{enumerate}
      Then, $F = \bigcup_{\delta > 0} \big( \XP \setminus \cup_{n = 1}^{\infty} A_{( n, \delta )} \big)$ is a cylinder absolute winning set with winning dimension at least $\beta_0$ where $\beta_0 := 1 / p_{i}$ if $i > 0$ and $1/3$ otherwise.
    \end{theorem}
    \begin{proof}
      Given any $0 < \beta < \beta_0$, let $B_0 = \overline{B}(\mathbf{x}, r_0 )$ be the init{i}al closed ball of radius $r_0$ chosen by Bob to kick start the cylinder $\beta$-absolute game. Without loss of generality, we may assume that $r_0 < 1/2$ as well as that the balls $B_i$ chosen by Bob have radii $r_i \rightarrow 0$ (Alice can force this by removing some largest possible cylinder which is legally allowed at each turn). We let $R,\ \eps$ and $( R_n )$ take the values dictated by our hypothesis for $K = B_0$ and $\mu = \beta^2/2$. Then, let $k_0 \in \N$ be the smallest such that $\mu^{k_0} < \min \{ \eps \mu r_0^{-1}, R^{-1} \}$ and $\delta := \mu^{k_0 + 1} r_0 < \eps$. For $k \geq 1$, mark $h_k \rightarrow \infty$ to be any strictly increasing subsequence such that
      \begin{equation}\label{E:h_k}
        \beta \mu^k r_0 < \mathrm{radius}\,( B_{h_k} ) =: r_k \leq \mu^k r_0.
      \end{equation}
      This is well-def{i}ned as $r_{k + 1} \geq \beta r_k$ for all $k \in \N$ and $\mu^{k + 1} < \beta \mu^k$. We claim that Alice is able to play in such a manner that the closed ball $B_{h_k}$ does not intersect any $A_{( n, \delta )}$ with $R_n \geq \mu^{k - k_0}$. The limit point $b_{\infty} = \cap_{k = 0}^{\infty} B_k = \cap_{k = 0}^{\infty} B_{h_k}$ shall then be in $F$ and the proof of the theorem will be done (as $\beta < \beta_0$ is arbitrary).\\[-0.2cm]

      Our claim is vacuously true for $k = 0$ as $R_n \leq R < \mu^{-k_0}$ for all $A_{( n, \delta )}$ intersect{i}ng $B_0$ non-trivially by our assumpt{i}on. Thereaf{t}er, supposing that the claim holds for $k$, we show it to be true for $k + 1$. Since the sets $A_{( n, \delta )}$ with the corresponding cylinder radii bounds $R_n \geq \mu^{k - k_0}$ have already been taken care of, we only need to show that Alice can now ensure $B_{h_{k + 1}}$ does not intersect $A_{( n, \delta )}$ for any $n \in N$ such that $\mu^{k + 1 - k_0} \leq R_n < \mu^{k - k_0}$. As hinted before, she has to worry about exactly one such subset. For, if both $A_{(n_1, \delta )} \cap B_{h_k},\ A_{( n_2, \delta )} \cap B_{h_k} \neq \phi$ and $R_{n_1}, R_{n_2} \in [\,\mu^{k + 1 - k_0},\,\mu^{k - k_0} )$, then the second condit{i}on of the theorem says that $d ( A_{( n_1, \delta )},\ A_{( n_2, \delta )} ) \geq \eps ( R_{n_1} + R_{n_2} ) \geq 2\eps \mu^{k + 1 - k_0}$ while $| B_{h_k} | \leq 2 r_k \leq 2\mu^k r_0$ and we have a contradict{i}on.\\[-0.2cm]
      
      If $n \in N$ is the unique index for which $A_{( n, \delta )} \cap B_{h_k} \neq \phi$ and $R_n \in [\,\mu^{k + 1 - k_0}, \mu^{k - k_0} )$ where $B_{h_k} = \overline{B} ( \x_k, r_k )$, Alice chooses $C_{h_k + 1}$ to be the open cylinder $C_{( n, \delta )}$ and since
      \begin{equation}\label{E:Aradius}
        \mathrm{radius}\,( C_{h_k + 1} ) \leq \delta R_n < \mu^{k_0 + 1} r_0\cdot\mu^{k - k_0} < \beta\cdot r_k,
      \end{equation}
      this const{i}tutes a legal move. It only remains to be argued that Bob has some choice of $B_{h_k + 1} \subset B_{h_k} \setminus C_{h_k + 1}$ lef{t} (in fact, plenty of them). If the constraining coordinate $i$ of $C_{h_k + 1}$ equals zero, we only need to find a point in the closed ball $\overline{B} \big( x_{k, 0}, ( 1 - \beta ) r_k \big) \subset \R$ which is at a Euclidean distance $\beta r_k$ from some open ball $B ( y, \beta r_k )$ containing the project{i}on $\pi_0 ( C_{h_k + 1} )$ of $C_{( n, \delta )}$ in the archimedean coordinate. This is clearly possible as long as $\beta < 1/3$.\\[-0.2cm]
      
      Else if $i > 0$, let $\overline{B} ( x_{k, i}, \lfloor p_i r_k \rfloor_i ),\ B ( y, p_i r'_{k + 1} ) \subset \Q_{p_i}$ be the respect{i}ve images of $B_{h_k}$ and $C_{h_k + 1}$ under the project{i}on $\pi_i$ onto the $i$-th coordinate. As $B_{h_k} \cap C_{h_k + 1} \neq \phi$ by our assumption, we get that $\overline{B} ( x_{k, i}, \lfloor p_i r_k \rfloor_i )\,\cap\,B ( y, r'_{k + 1} ) \neq \phi$ too. Being balls in an ultrametric space, one of them then has to be contained in the other and because $\beta < 1 / p_i$, we have
      \begin{equation}\label{E:ultra}
        r'_{k + 1} \leq \frac{1}{p_i} \lfloor r_k \rfloor_i
      \end{equation}
      and thereby $B ( y, p_i r'_{k + 1} ) \subsetneq \overline{B} ( x_{k, i}, p_i \lfloor r_k \rfloor_i ) = \overline{B} ( y, p_i \lfloor r_k \rfloor_i )$. Bob picks a point $z_i \in \overline{B} ( y, p_i \lfloor r_k \rfloor_i )$ whose distance from $y$ is equal to $p_i \lfloor r_k \rfloor_i$ and if
      \begin{equation}\label{E:cproj}
        ( x_{k, 0}, \ldots, x_{k, i - 1}, x_{k, i + 1}, \ldots ) = \pi^{\perp}_i ( \mathbf{x}_k ),
      \end{equation}
      where $\pi^{\perp}_i : X \rightarrow \R \times \prod^{\prime}_{j \neq i} \Q_{p_j}$ is the complementary project{i}on of $\pi_i$, he def{i}nes
      \begin{equation}\label{E:nextball}
        B_{h_k + 1} = \overline{B} \big( ( x_{k, 0}, \ldots x_{k, i - 1}, z_i, x_{k, i + 1}, \ldots ), \beta r_k \big).
      \end{equation}
      This has diameter $2 \beta r_k$, is contained in $B_{h_k}$ and avoids the cylinder $C_{h_k + 1}$.
    \end{proof}
    {Given any posit{i}ve lower bound on the \CAW\ dimension, we can boost it to an absolute quant{i}ty (depending on \P\ alone) for f{i}nite solenoids.
    \begin{lemma}\label{L:windim}
      Let $\P = \{ p_1 < \cdots < p_{l - 1} \}$ be f{i}nite. Any \CAW\ subset $S$ of \XP\ has winning dimension $\geq \min \betaP := \{ 1/3, 1 / p_{l - 1} \}$.
    \end{lemma}
    \begin{proof}
      Assume $0 < \beta < \betaP$. We are guaranteed that there exists some $0 < \beta' < \beta$ such that $S$ is cylinder $\beta'$-absolute winning. Changing the game parameter from $\beta'$ to $\beta$ only enlarges the set of choices available to Alice while Bob cont{i}nues to have some legal choice lef{t} as long as $\beta < 1/3$ and $\beta < 1 / p_{l - 1} \leq 1 / p_i$, if $i > 0$ is the constraining coordinate of the cylinder blocked by Alice in the previous move. Also, all of his valid moves in the cylinder $\beta$-absolute game remain so in the $\beta'$-game. Alice just needs to pretend that the game parameter is $\beta'$ and follow her winning strategy for the same.
    \end{proof}}
    
  \section{Non-dense orbits of solenoidal maps}\label{S:morph}
    As already men{t}{i}oned in \S\,\ref{S:intro}, {an af{f}ine endomorphism $A: \SP \rightarrow \SP$ is of the form $\x \mapsto (m/n)\x + \a$ for some $m/n \in R$ and $\a \in \SP$.} Here, $R$ is the set of endomorphisms of the solenoid $\SP = \nicefrac{\XP}{\Delta(R)}$ given by the ring $R = \Z \big[ \{ 1/ p_i \mid p_i \in \P \} \big]$. {The af{f}ine transformat{i}on $A$ is invert{i}ble if{f} $n/m \in R$ too.}\\[-0.1cm]
    
    Next, the cylinder $\beta$-absolute game on \SP\ can be shi{f}{t}ed to a game played on \XP\ once the radii of the balls $B_i \subset \SP$ become small enough (say $< 1/2$). This can be forced on Bob in {f}{i}nitely many steps after the beginning of the game.\\[-0.1cm]

    Pick some $y \in \SP$ and let $A$ be any {f}{i}xed af{f}ine transformat{i}on of \SP\ with its linear part $m/n \in R \setminus \{ 0, \pm 1 \}$. We abuse nota{t}{i}on and call any of its li{f}{t}s from $\XP \rightarrow \XP$ to be $A$ too. Note that any such lif{t} is an invert{i}ble self map of \XP\ as long as $A: \SP \rightarrow \SP$ is surject{i}ve. Further, let $F_A (y)$ denote the set of points $\mathbf{x} \in \XP$ whose image $x := \Pi ( \mathbf{x} )$ has its $A$-orbit not entering some $\delta ( x, y, A )$-neighbourhood of $y$. The goal for Alice is to avoid $\eps$-neighbourhoods of the grid points $\Pi^{-1} ( \{ y \} ) = \mathbf{y} + \Delta(R)$ (for some $\mathbf{y} \in \Pi^{-1} ( \{ y \} )$) which are all at least a unit distance away from each other. Otherwise said, the set $F_{A} (y)$ that Alice should aim for is
    \begin{equation}\label{E:winset}
      \bigcup_{t > 0} \left( \XP \setminus \bigcup_{j = 0}^{\infty} A^{-j} \big( \Delta (R) + B ( \mathbf{y}, t ) \big) \right)
    \end{equation}
    and $\nicefrac{F_{A} (y)}{\Delta(R)}$ shall be the image set for the game played on \SP. By our assump{t}{i}on about $A$, there exists $i \geq 0$ such that $| \nicefrac{m}{n} |_{p_i} =: \lambda_i > 1$. We let $\lambda_A$ to be $\sup_i \lambda_i$. This is f{i}nite, at{t}ained for some $i = i_0$ and strictly greater than $1$. In part{i}cular for any $\mathbf{x}_1, \mathbf{x}_2 \in \XP$, we have that
    \begin{equation}\label{E:dlb}
      d (\, A^{-j} \mathbf{x}_1,\, A^{-j} \mathbf{x}_2 \,) \geq \lambda_A^{-j} d (\, \mathbf{x}_1,\, \mathbf{x}_2 \,) \textrm{ for all } j \geq 0
    \end{equation}
    {as translat{i}on by any \a\ is an isometry of \SP\ (and \XP). Equivalently for any two subsets $F_1, F_2 \subset \XP$,
    \begin{equation}
      d (\, A^{-j} F_1,\, A^{-j} F_2 \,) \geq \lambda_A^{-j} d (\, F_1,\, F_2 \,)\ \forall j \geq 0.
    \end{equation}}
    The constraining coordinate of all the cylinders removed by Alice will be some f{i}xed $i_0$ for which $\lambda_{i_0} = \lambda_A$. Let $0 < \mu < 1$ and $\ell \in \N$ be the smallest for which $\lambda_A^{-\ell} < \mu$. If $a = m^{-\ell}$, then $(\nicefrac{m}{n})^{-j}R \subseteq aR$ for all $j \in \{ 0, \ldots, \ell \}$. Not unlike $R$, the points of $aR$ too const{i}tute a $\delta$-uniformly discrete set for some $0 < \delta = \delta (a, \P) \leq 1$. We also choose $b \geq 1$ given by
    \begin{equation}
      b = \max_{0 \leq j \leq \ell} \sup_{\x \in \overline{B ( \mathbf{0}, 1 )}} d \big( A^{-j} \x, \mathbf{0} \big)
    \end{equation}
    and let
    \begin{equation}\label{E:t0}
      t_0 = \min \dfrac{1}{3b} \big\{ d (\,\mathbf{y} - A^{-j}\mathbf{y},\,a\mathbf{z}\,) > 0 \mid 0 \leq j \leq \ell,\ \mathbf{z} \in \Delta(R) \big\}.
    \end{equation}
    This belongs to $]\,0,\,1/3\,]$ and thereby for any $\mathbf{z}_1, \mathbf{z}_2 \in \Delta(R)$ such that $\mathbf{y} + \mathbf{z}_1 \neq A^{-j} ( \mathbf{y} + \mathbf{z}_2 )$ for some $0 \leq j \leq \ell$, we have
    \begin{align}\label{E:dcalc}
      d \left( B( \mathbf{y} + \mathbf{z}_1, t_0 ), A^{-j} ( B( \mathbf{y} + \mathbf{z}_2, t_0 ) ) \right) &\geq d \left( B( \mathbf{y} + \mathbf{z}_1, bt_0 ), B( A^{-j} (\mathbf{y} + \mathbf{z}_2), bt_0 ) \right)\notag\\
              &\geq 3bt_0 - 2bt_0 = bt_0 \geq t_0.
    \end{align}
    If $j_1 \leq j_2$ are any two exponents such that $\mu^{k + 1} < \lambda_A^{-j_2} \leq \lambda_A^{-j_1} \leq \mu^k$ for some $k \in \N$, then $j_2 - j_1 \leq \ell$ by the very def{i}nit{i}on of $\ell$. Hence, for $0 < t < \mu t_0 / 2$ and any $\mathbf{z}_1, \mathbf{z}_2 \in \Delta(R)$ for which $\mathbf{y} + \mathbf{z}_1 \neq A^{- ( j_2 - j_1 )} ( \mathbf{y} + \mathbf{z}_2 )$, we get that
    \begin{equation}\label{E:farset}
      d \left( A^{-j_1} B ( \mathbf{y} + \mathbf{z}_1, t ), A^{-j_2} ( B ( \mathbf{y} + \mathbf{z}_2, t ) ) \right) > \mu^{k + 1} t_0 = \frac{\mu t_0}{2} ( \mu^k + \mu^k )
    \end{equation}
    while $A^{-j} B ( \mathbf{y} + \mathbf{z}, t ) \subset C \big( A^{-j} ( \mathbf{y} + \mathbf{z} ), \lambda_A^{-j}t, i_0 \big) \subseteq C \big( A^{-j} ( \mathbf{y} + \mathbf{z} ), \mu^k t, i_0 \big)$ for $\lambda_A^{-j} \leq \mu^k$. We let
    \begin{equation}\label{E:index}
      N = \{\,(\,j,\,\mathbf{y} + \mathbf{z}\,) \mid j \in \N,\ \mathbf{z} \in aR\,\}
    \end{equation}
    be the countable indexing set in our Theorem~\ref{Th:cgame}. The f{i}rst hypothesis therein is satisf{i}ed by taking $R = 1$ and for $n = ( j, \mathbf{y} + \mathbf{z} )$, let{t}ing
    \begin{equation}
      A_{( n, t )} = A^{-j} B ( \mathbf{y} + \mathbf{z}, t )\ \subset\ C_{( n, t )} :=  C \big( A^{-j} ( \mathbf{y} + \mathbf{z} ), \lambda_A^{-j}t, i_0 \big)
    \end{equation}
    which suggests that we should take $R_n = \lambda_A^{-j}$. Clearly, the second hypothesis has been shown to hold here in~\eqref{E:farset}. Hence, we infer that $F_{A} (y) = \cup_{t > 0} \big( \XP \setminus \bigcup_{j = 0}^{\infty} A^{-j} \big( R + B ( \mathbf{y}, t ) \big) \big)$ is a \CAW\ subset of \XP\ with winning dimension as in the statement of Theorem~\ref{Th:cgame} and so is its image in \SP. Because Proposit{i}on~\ref{P:ccwin}, we can extend this result to the set of points whose $A$-orbits avoid some neighbourhoods of countably many points $\{ \y_k \}_{k \in \N} \subset \SP$.\\[-0.1cm]
    
    If $A (\x) = (m/n)\x + \a$ is such that $m/n = -1$, then $A^2$ is the ident{i}ty endomorphism. In this case, Alice only needs to move the game away from the countable set $\{ \y_k \} \cup \{ \a - \y_k \}$. This is trivial. The situat{i}on is even simpler when $A$ is just the ident{i}ty map. Now, let $Y$ be the set consist{i}ng of all those points of \SP\ which have a periodic orbit for some $B \in R \setminus \{ \pm 1 \}$. This is countable and leads us to conclude:
    \begin{theorem}\label{Th:solmain}
      {Let $\{ A_j : \x \mapsto ( m_j / n_j ) \x + \a_j \}_{j \in \N}$ be any subset of af{f}ine surject{i}ve endomorphisms of the solenoid \SP\ such that
      \begin{enumerate}
        \item none of the $A_j$'s is a non-trivial translat{i}on of \SP, and
        \item the collect{i}on of rat{i}onal numbers $\{ m_j / n_j \}_{j \in \N}$ belong to some f{i}nite extension of \Z.
      \end{enumerate}}
      Then, the set of points whose orbit closure under the act{i}on of any of the $A_j$'s does not contain any periodic $B$-orbit for all $B \in R \setminus \{ \pm 1 \}$ is cylinder absolute winning.
    \end{theorem}
    \begin{proof}
      If $\{ A_j \} \subset \Z \big[ \{ 1 / p_1, \ldots, 1 / p_n \mid p_i \in \P \} \big]$, then the winning dimension of each of the subsets $F_{A_j}$ is at least $\min \{ 1/3, \min \{ 1 / p_i \mid 1 \leq i \leq n \} \} > 0$. This is also a lower bound on the \CAW\ dimension of the intersect{i}on $\cap_j F_{A_j}$ invoking Proposit{i}on~\ref{P:ccwin} once again.
    \end{proof}
    Note that even though $R$ is a countable set, we cannot further this argument to take intersect{i}ons over any arbitrarily chosen sequences of af{f}ine surject{i}ve endomorphisms of \SP. This is because the lower bound on the winning dimension of the \CAW\ subsets of \SP\ corresponding to each $A$ is dependent on $A$ itself in terms of $i_0$ for which $\lambda_{i_0} = \lambda_A$. However, for f{i}nite \P, each such $\beta_0$ is at least $\min \{ 1/3, \min \{ 1 / p_i \mid p_i \in \P \} \}$. We can then remove the second condit{i}on in Theorem~\ref{Th:solmain} to get
    \begin{theorem}\label{Th:sol}
      Let \P\ be a f{i}nite set of rat{i}onal primes {and $\{ A_j \}$ be any sequence of af{f}ine surject{i}ve endomorphisms of \SP\ such that none of the $A_j$'s is a translat{i}on.} The set of points whose orbit closure under the act{i}on of any $A_j$ does not contain any periodic $B$-orbit for all $B \in R \setminus \{ \pm 1 \}$ is \CAW\ with winning dimension at least $\min \{ 1/3, 1 / p_{l - 1} \}$. Here, $p_{l - 1}$ is the largest prime in \P.
    \end{theorem}
    {\noindent In part{i}cular, this is true of the collect{i}on of all surject{i}ve endomorphisms of \SP.}
    
  \section{\protect{Sizes of \CAW\ subsets}}\label{S:hausdorff}
    {Let \P\ be f{i}nite. We start by discussing the implicat{i}ons of \CAW\ property of a subset $F$ for a strong game played on \XP\ with $F$ as its target.
    \begin{proposition}\label{P:CAWstrong}
      A \CAW\ subset of \XP\ is $\alpha$-strong winning for all $\alpha < \betaP$.
    \end{proposition}}
    \begin{proof}
      Without loss of generality, we may take the \CAW\ dimension of $F$ to be $\betaP$ due to Lemma~\ref{L:windim}. This means our target set $F$ is $\beta$-CAW for all $0 < \beta < \betaP$. Now, suppose that $\alpha \in\ ]\,0,\,\betaP\,[$ and $\gamma \in\ ]\,0,\,1\,[$ are any f{i}xed (strong) game parameters for Alice and Bob, respect{i}vely.\\[-0.2cm]
      
      Given a ball $B_0 = \overline{B} ( \mathbf{x}, r )$ chosen by Bob at any stage of the strong game, Alice checks the cylinder $C$ with $\Radius{C} \leq \alpha\gamma r$ to be removed by her in accordance with her winning strategy for $F$ when playing the cylinder $( \alpha\gamma )$-absolute game. If $B \cap C = \phi$, she chooses any $A \subset B$ allowed by the rules of the strong game. Assume this to not be the case for the rest of this proof.\\[-0.2cm]
      
      If the constraining coordinate $i$ of $C$ is archimedean, Alice has no problem in choosing a Euclidean ball $A \subset \pi_0 (B) \setminus \pi_0 (C)$ with $\Radius{A} \geq \alpha r$ as $\alpha \leq \betaP \leq 1/3$. Else again when $i > 0$, we have $\pi_i ( C ) \subsetneq \pi_i (B)$ because
      \begin{equation}
        \Radius{\pi_i (C)} \leq p_i \floor{\betaP r}_i \leq \floor{r}_i \text{ while } \Radius{\pi_i (B)} \geq p_i \floor{r}_i.
      \end{equation}
      We can moreover take the center $x_i$ of $\pi_i (B)$ to be the same as that of $\pi_i (C)$. Let $z_i \in \pi_i (B) \setminus \pi_i (C)$. Then, $\abs{z_i - x_i }_{p_i} = p_i \floor{r}_i$ and the ultrametric also gives us that $\overline{B}_{\Q_{p_i}} ( z_i, \floor{r}_i ) \subset \pi_i (B) \setminus \pi_i (C)$. In either case, the pre-images $\pi_0^{-1} (A)$ or $\pi_i^{-1} \big( \overline{B}_{\Q_{p_i}} ( z_i, \floor{r}_i ) \big)$ contain a ball of \XP\ of radius at least $\betaP r \geq \alpha r$ which lies inside $B \setminus C$. Alice chooses one such $A_1$ to be her next move. Bob's choice of any $B_1 \subset A_1$ with $\Radius{B_1} \geq \gamma\cdot\Radius{A_1} \geq \alpha\gamma r$ immediately af{t}er is also a valid move in cylinder $(\alpha\gamma)$-absolute game.
    \end{proof}
    It should be ment{i}oned here that the relat{i}onship between winning sets for strong games and quasisymmetric homeomorphisms of \XP\ is not clear to us. Nor do we have the analogous statement of Proposit{i}on~\ref{P:CAWstrong} for the full solenoid.
    
    \subsection{Incompressibility}\label{SS:incom}
    The next result is about the incompressible behaviour of cylinder absolute winning subsets of \XP. A set $S \subset \XP$ is \emph{strongly af{f}inely incompressible} if for any non-empty open subset $U$ and any sequence of invert{i}ble af{f}ine homomoprhisms $( \Psi_i )_{i \in \N}$, the set $\cap_{i \in \N} \Psi_i^{-1} S \cap U$ has the same Hausdorf{f} dimension as $U$~\citep{BFK+12, Dan89}. It is our claim that \CAW\ subsets of \XP\ are strongly af{f}inely incompressible for f{i}nite \P. We show this by proving a lower bound on the \CAW\ dimension of $\Psi^{-1} S \cap U$ in terms of $\operatorname{windim} S$ for any af{f}ine map $\Psi$ of \XP. Together with Proposit{i}on~\ref{P:ccwin}, this will give us that the intersect{i}on of any countably many pre-images of a \CAW\ subset under invert{i}ble af{f}ine homomorphisms is \CAW\ too.
    \begin{theorem}
      Let \P\ be f{i}nite, $U \subset \XP$ open and $S$ be any \CAW\ subset. Also, let $\Psi : \XP \rightarrow \XP$ be an invert{i}ble af{f}ine homomorphism. Then, the set $\Psi^{-1} S \cup ( \XP \setminus U )$ is also \CAW\ with winning dimension at least $\betaP$.
    \end{theorem}
    \begin{proof}
      As in Proposit{i}on~\ref{P:CAWstrong}, we may take the \CAW\ dimension of $S$ to be $\betaP$ without loss of generality. Let us f{i}rst make some reduct{i}ons to simpler situat{i}ons. If the diameters $\abs{B_k}$ of balls chosen by Bob don't go to zero as $k \rightarrow \infty$, then $\cap_k B_k$ contains an open ball inside it. As $S$ is a winning subset, it has to be dense and in turn its pre-image $\Psi^{-1} S$ is also dense in \XP. Second, if $B_k \cap ( \XP \setminus U ) \neq \phi$ for inf{i}nitely many $k$, then they form a decreasing sequence of closed subsets of the compact ball $B_1$. Their intersect{i}on cannot be empty and hence $\cap_k B_k$ contains a point of $\XP \setminus U$ result{i}ng in Alice's victory. It is safe to exclude both of these events from the rest of the proof. We can moreover take that $B_0 \subset U$.\\[-0.2cm]
      
      Let $0 < \beta < \betaP$ and $\Psi (\x) = D\x + \a$ as explained before. Following~\citep{BFK+12}, Alice will run a `hypothet{i}cal' Game~2 (in her mind) where the target set is $S$ and a dif{f}erent game parameter $\beta'$ which is some posit{i}ve power of $\beta$. She carefully decides and projects some of Bob's moves in the $\big( \Psi^{-1} (S) \cup ( \XP \setminus U ), \beta \big)$-game to construct choices made by a hypothetical Bob~II in Game~2. Since there is a winning strategy for the lat{t}er by hypothesis, she channels the winning moves in this second game via the inverse map $\Psi^{-1}$ to win over $\Psi^{-1} (S)$. We take
      \begin{equation}
        \lambda_{\Psi} := \max \{\,\max_{p \in \P} \abs{D}_p,\,\abs{D}\,\}
      \end{equation}
      which makes sense as $D$ is a rat{i}onal number. As $D \neq 0$, we have $1 \leq \lambda_{\Psi} < \infty$ for any $\Psi$ and any \P. It is clear that
      \begin{equation}
        d \big( \Psi (\x), \Psi (\y) \big) \leq \lambda_{\Psi} d ( \x, \y )
      \end{equation}
      as $d$ is a translat{i}on-invariant metric on \XP. We re-label the choices made by Bob such that $\Radius{B_0} < 1 / \lambda_{\Psi}$. Let $n \in \N$ be the smallest posit{i}ve natural number for which
      \begin{equation}
        \lambda_{\Psi}\lambda_{\Psi^{-1}} ( \beta + 1 ) \beta^{n - 2} < 1,\quad \beta' = \beta^n\quad \textrm{and}\quad \eta := ( \beta + 1 ) \beta^{n - 1}.
      \end{equation}
      Alice waits for the stages $0 = j_1 < j_2 < \cdots$ in the original $\big( \Psi^{-1} (S) \cup ( \XP \setminus U ), \beta \big)$-cylinder absolute game when
      \begin{equation}
        \beta^{n} \leq \Radius{B_{j_k}} / \Radius{B_{j_{k - 1}}} < \beta^{n - 1}
      \end{equation}
      for the f{i}rst t{i}me. Not{i}ce that this is well-def{i}ned and exists because we assumed $\abs{B_k} \rightarrow 0$ and the radius of Bob's choice at $( k + 1 )$-th step cannot shrink by a factor of more than $\beta$ compared to that of his choice at the $k$-th step for all $k$. Imitat{i}ng~\cite{BFK+12}, denote
      \begin{equation}
        B_{j_k} = \overline{B} ( \x_k, r_k ) \textrm{ and } B'_k = \overline{B} ( \Psi (\x_k), r'_k ) \textrm{ where } r'_k = \lambda_{\Psi} r_k
      \end{equation}
      so that $\Psi ( B_{j_k} ) \subset B'_k$ for all $k$ by our def{i}nit{i}on of $\lambda_{\Psi}$. Then, $\Psi ( \cap_k B_k ) = \Psi ( \cap_k B_{j_k} ) \subset \cap_k B'_k$ while the intersect{i}ons $\cap_k B_k$ and $\cap_k B'_k$ are both singleton sets (we are in the case when $\abs{B_k} \rightarrow 0 \Rightarrow \abs{B'_k} \rightarrow 0$ as $k \rightarrow \infty$). Thus, we see that the non-empt{i}ness of $\cap_k B'_k \cap S$ will imply that of $\cap_k B_k \cap \Psi^{-1} (S)$.\\[-0.2cm]
      
      If $C'_{k + 1} = C ( \y'_{k + 1}, \eps'_{k + 1}, i_{k + 1} )$ is the cylinder to be removed by Alice in Game~2 where $\eps'_{k + 1} \leq \beta' r'_{k}$, she chooses $C_{j_k + 1}$ as
      \begin{equation}
        C \left( \Psi^{-1} ( \y'_{k + 1} ), \lambda_{\Psi}\lambda_{\Psi^{-1}}\eta r_{k}, i_{k + 1} \right) \supset \Psi^{-1} \big( C ( \y'_{k + 1}, \eps'_{k + 1}, i_{k + 1} ) \big)
      \end{equation}
      to be blocked next in the $\beta$-cylinder absolute game with target set $\Psi^{-1} (S)$. By design, $\lambda_{\Psi}\lambda_{\Psi^{-1}}\eta < \beta$ and it only remains to show that Bob has some choice of $B_{j_k + 1} \subset B_{j_k} \setminus C_{j_k + 1}$ lef{t} with $\Radius{B_{j_k + 1}} \geq \beta\cdot\Radius{B_{j_k}}$. As we have seen in the proof of Theorem~\ref{Th:cgame}, this is clearly not a problem as $\beta < \betaP \leq \min \{ 1/3, 1 / p_{i_{k + 1}} \}$ when $i_{k + 1} > 0$ or otherwise. All of Bob's subsequent choices in Game~1, including $B_{j_{k + 1}} = \overline{B} ( \x_{k + 1}, r_{k + 1} )$, obey
      \begin{equation}
        \abs{x_{k, 0} - x_{k + 1, 0}} \leq r_k - r_{k + 1}
      \end{equation}
      in the archimedean coordinate and
      \begin{equation}
        \abs{x_{k, i} - x_{k + 1, i}}_{p_i} \leq p_i r_k\ \forall i > 0.
      \end{equation}
      Then,
      \begin{align}
        \abs{\Psi ( \x_{k + 1} )_0 - \Psi ( \x_{k} )_0} &= \abs{Dx_{k + 1, 0} + a_0 - (Dx_{k, 0} + a_0)}\\
        &\leq \lambda_{\Psi} \abs{x_{k + 1, 0} - x_{k, 0}} \leq \lambda_{\Psi} ( r_k - r_{k + 1} ) = r'_k - r'_{k + 1},\notag
      \end{align}
      and for all $i > 0,\ \abs{\Psi ( \x_{k + 1} )_i - \Psi ( \x_{k} )_i}_{p_i} \leq \lambda_{\Psi} p_i r_k = p_i r'_k$ similarly. The conclusion cannot be escaped that the corresponding ball $B'_{k + 1} = \overline{B} \big( \Psi(\x_{k + 1}), r'_{k + 1} \big) \subset B'_k$ in Game~2. It is also outside of $C'_{k + 1}$ when $i_{k + 1} > 0$ as
      \begin{equation}
        \abs{x_{k + 1, i_{k + 1}} - \Psi^{-1} ( \y'_{k + 1} )_{i_{k + 1}}}_{i_{k + 1}} \geq p_{i_{k + 1}} \lambda_{\Psi}\lambda_{\Psi^{-1}}\eta r_k
      \end{equation}
      which gives that
      \begin{align}
        \abs{\Psi ( \x_{k + 1} )_{i_{k + 1}} - y'_{k + 1, i_{k + 1}}}_{i_{k + 1}} &\geq p_{i_{k + 1}} \lambda_{\Psi}\eta r_k = p_{i_{k + 1}} \eta r'_k\notag\\
        &= p_{i_{k + 1}} ( \beta' + \beta^{n - 1} ) r'_k\\
        &\geq p_{i_{k + 1}} ( \eps'_{k + 1} + r'_{k + 1} )\notag
      \end{align}
      by Alice's choice of marking for $B_{j_{k + 1}}$. The computat{i}ons are not very dif{f}erent when the constraining coordinate of $C'_{k + 1}$ is archimedean.
    \end{proof}
    For general \P, we are only able to show the largeness of countable intersect{i}ons of pre-images under translations of \XP.
    \begin{proposition}\label{P:translation}
      Let $( \a_k )_{k \in \N} \subset \XP$ be any arbitrary sequence and $S$ be a \CAW\ subset with winning dimension $\beta_0$. Then, so is $S \cap \bigcap_{k \in \N} ( S + \a_k )$.
    \end{proposition}
    \begin{proof}
      Each of the translat{i}ons $\Psi_k (\x) := \x - \a_k$ is an isometry and in part{i}cular, does not change the shape of the balls in \XP. Given any such single $\Psi$, we argue that $\Psi^{-1} (S)$ has the same \CAW\ dimension as $S$. Alice simply translates back her choices for the Game~2 described above by $-\a$ when $\Psi (\x) = \x + \a$ and projects Bob's succeeding choice forward by $\Psi$. Note that as $\lambda_{\Psi_k} = \lambda_{\Psi_k^{-1}} = 1$ for all $k$, she should take $\beta' = \beta$ for any $0 < \beta < \beta_0$. One should also replace $\eta = 1$ in the previous calculat{i}ons. The countable intersect{i}on property then follows by Proposit{i}on~\ref{P:ccwin}.
    \end{proof}
    
    \subsection{\protect{Hausdorf{f} dimension and measures}}\label{SS:Hau}
      Lastly, we will try to understand the sizes of \CAW\ sets in terms of Hausdorf{f} dimensions and measures. Towards this goal, we will require an est{i}mate on the number of legal choices that Bob has at any stage of the game.
      \begin{lemma}\label{L:nodisjointballs}
        Let $0 < \beta \ll 1$. Then, the maximum number of pairwise disjoint balls of radius $\beta r$ contained in any closed ball $\overline{B} ( \mathbf{x}, r ) \subset \XP = \R \times \prod^{\prime}_{j > 0} \Q_{p_j}$ which do not intersect an open cylinder $C ( \mathbf{y}, \beta r, i )$ and also maintain a distance at least $\beta r$ from each other is given by
        \[
          N_{\C} ( \beta ) \gg \beta^{- \frac{\thetaP{1 / \beta}}{\ln ( 1 / \beta )}}
        \]
        where $\thetaP{t} = \sum_{p \in \P,\ p \leq t} \ln p$.
      \end{lemma}
      \begin{proof}
        Because~\eqref{E:solmet}, every closed ball is the Cartesian product of its coordinatewise projections. In any non-archimedean coordinate $j$ such that $p_j \leq 1 / \beta$, there are at least $( p_j \floor{\beta}_j )^{-1}$ pairwise disjoint balls contained in the project{i}on $\pi_j \big( \overline{B} ( \mathbf{x}, r ) \big) = \overline{B} ( x_j, \lfloor p_j r \rfloor_j ) $ whose radius equals $\lfloor p_j \beta r \rfloor_j$. Each of them also maintain a distance of at least $p_j^2 \floor{\beta r}_j$ from each other which means that the pre-images of any two such sub-balls in \XP\ are $\geq p_j \floor{\beta  r}_j > \beta r$ away. The lower bound $( p_j \floor{\beta}_j )^{-1}$ equals $\floor{1 / \beta}_j$ unless $\beta$ is an integral power of $p_j$ in which case it is $p_j^{-1} \floor{1 / \beta}_j$. Note that this can happen for at most one prime for any given $\beta$ and we have already assumed $\beta < 1 / p_j$. In the real coordinate, this number is $\gg \beta^{-1}$ even when we ask that the balls are at least $\beta r$ apart. The pre-images under the project{i}on map $\pi_{\beta} : \XP \rightarrow \R \times \prod_{p_j \leq 1 / \beta} \Q_{p_j}$ of any product of these sub-balls of $\pi_j \big( \overline{B} ( \mathbf{x}, r ) \big)$ for $j = 0$ or $p_j \leq 1/\beta$ are pairwise disjoint, each contain at least one sub-ball of $\overline{B} ( \mathbf{x}, r )$ of radius $\beta r$ and the minimum distance between any two of those pre-images is $\geq \beta r$.\\[-0.2cm]
      
      When asking for only those sub-balls that do not intersect the open cylinder $C ( \mathbf{y}, \beta r, i )$, it is necessary and suf{f}icient that we restrict ourselves to only those from the above chosen collect{i}on whose images in $\pi_i \big( QB^{-} ( \mathbf{x}, r ) \big)$ do not intersect $\pi_i \big( C ( \mathbf{y}, \beta r, i ) \big)$. Otherwise said, all coordinates but $i$ are not af{f}ected. If $i = 0$, the number of such balls in $\pi_0 \big( QB^- ( \mathbf{x}, r ) \big) = \overline{B} ( x_0, r )$ that do not intersect $\pi_0 \big( C ( \mathbf{y}, \beta r, 0 ) \big) = B ( y_0, \beta r )$ is still $\gg \beta^{-1}$ albeit with a smaller constant. Else, it is at least $( p_i \beta )^{-1} - 1 > \beta^{-1} / 2p_i$. The Cartesian product of these disjoint balls in $\Q_{p_j}$'s and \R\ then gives us that
      \begin{equation}\label{E:NCbeta}
        N_{\C} ( \beta ) \gg \prod_{\substack{p_j \in \P,\\ p_j < 1 / \beta}} \floor{\beta^{-1}}_j \geq \beta^{- \frac{\thetaP{1 / \beta}}{\ln ( 1 / \beta )}}\ \forall 0 < \beta \ll 1
      \end{equation}
      by a calculat{i}on similar to the one in Lemma~\ref{L:prod}.
    \end{proof}
    As discussed in \S\,\ref{S:sol}, $\thetaP{t} \sim t$ when $t \rightarrow \infty$ and \P\ is the full set of primes. More generally, it shows an asymptot{i}c linear growth with $t$ when \P\ is any (inf{i}nite) set of all primes in an arithmet{i}c progression. We record that the constant implied by the Vinogradov notat{i}on in~\eqref{E:NCbeta} is independent of the constraining coordinate of the cylinder $C$. For f{i}nite \P, we use a slightly dif{f}erent lower bound.
    \begin{lemma}\label{L:ndbfin}
      Let $\abs{\P} = l - 1$ and $p_{l - 1}$ be the largest prime in \P. Then,
      \[
        N_{\C} (\beta) \gg \beta^{-l}\ \forall 0 < \beta \ll 1 / p_{l - 1},
      \]
      where the implied constant may depend on the primes $p_1, \ldots, p_{l - 1}$ and $l$.
    \end{lemma}
    \begin{proof}
      We only need to replace the lower bound for the number of disjoint sub-balls in each non-archimedean coordinate by $\beta^{-1} / p_j$ and the rest of the argument remains the same.
    \end{proof}
    Suppose $F$ is a cylinder absolute winning subset of $X$ and that Alice always plays according to a winning strategy if it is available. We shall now construct a subset $F^* \subseteq F$ which corresponds to the points obtained when Bob is only to allowed to choose one of the $N_{\C} (\beta)$-many sub-balls described in Lemmata~\ref{L:nodisjointballs} or \ref{L:ndbfin} at each stage of the game. This resembles closely a device from \citeauthor{Kri06}~\cite{Kri06} (see also~\cite[Theorem~6]{Sch66}).
    \begin{proposition}\label{P:hdim}
      Let $f$ be any dimension funct{i}on such that
      \[
        \limsup_{\delta \rightarrow 0} \frac{\log f (\delta)}{\log \delta} < \liminf_{\beta \rightarrow 0} \frac{\log N_{\C} (\beta)}{\abs{\log \beta}}.
      \]
      Then, the $f$-dimensional Hausdorf{f} measure of any \CAW\ subset of \XP\ is greater than zero.
    \end{proposition}
    \begin{proof}
      Let $\beta_0, \delta_0 > 0$ be small enough so that $\beta_0$ is less than the winning dimension of our \CAW\ set $F$, $\delta_0 < 1$ and
      \begin{equation}
        \sup_{\delta < \delta_0} \big( \log f (\delta) / \log \delta \big) < \inf_{\beta \leq \beta_0} \big( \log N_{\C} (\beta) / \abs{\log \beta} \big).
      \end{equation}
      They exist by virtue of our hypothesis about $f$. Now, let $\Lambda := \{ 0, 1, \ldots, N_{\C} ( \beta_0 ) - 1 \}^{\N}$, the sequence space each of whose element $\lambda = ( \lambda_k )_{k \in \N}$ corresponds to a sequence of choices made by Bob when he is only allowed to choose from one of the $N_{\C} ( \beta_0 )$-many disjoint sub-balls inside $B_{k - 1}$. The choices made by him at the $k$-th stage are labelled $B ( \lambda_1, \lambda_2, \ldots, \lambda_k )$. If $\lambda \neq \lambda'$, they di{f}{f}er in some entry $k_0$ and the corresponding balls $B ( \lambda_1, \ldots, \lambda_{k_0} )$ and $B ( \lambda'_1, \ldots, \lambda'_{k_0} )$ are disjoint. This implies that the points obtained at inf{i}nity, $\mathbf{a}_{\infty} ( \lambda ) = \cap_{k \rightarrow \infty} B ( \lambda_1, \ldots, \lambda_k ) \neq \mathbf{a}_{\infty} ( \lambda' ) = \cap_{k \rightarrow \infty} B ( \lambda'_1, \ldots, \lambda'_k ) \in F$ under the belief that Alice is following a winning strategy for the target set $F$ with game parameter $\beta_0$. Let
      \begin{equation}\label{E:Fstar}
        F^* := \{ \mathbf{a}_{\infty} ( \lambda ) \mid \lambda \in \Lambda \} \subseteq F.
      \end{equation}
      {We show that $\H^f (F^*) > 0$} and this shall in turn give us our desired statement. For this, the space $F^*$ is mapped in a cont{i}nuous fashion (via the biject{i}on with $\Lambda$) onto $[ 0, 1 ]$ using the $N_{\C} ( \beta_0 )$-adic expansion of real numbers, namely $\mathbf{a}_{\infty} (\lambda) \mapsto 0.\lambda_1 \lambda_2 \cdots$. Call this map $\psi$ and let $( U_n )_{n \in \N}$ be any $\delta$-cover of $F^*$ for some $\delta < \delta_0$. Without loss of generality, let $U_n \subset F^*$ for all $n$. Plainly, $\big( \psi ( U_n ) \big)_{n \in \N}$ is a cover for $[ 0, 1 ]$ and since diameter is an outer measure on \R, we get that
      \begin{equation}
        1 \leq \sum_{n \in \N} \abs{\psi ( U_n )}.
      \end{equation}
      Def{i}ne
      \begin{equation}
        j_n = \left\lfloor \dfrac{\log ( 2 \abs{U_n} )}{\log \beta_0} \right\rfloor
      \end{equation}
      so that $j_n > 0$ for all $\abs{U_n}$ small enough and furthermore, $\abs{U_n} < \beta_0^{j_n}$. Thus, $U_n$ intersects non-trivially with at most one of the balls $B ( \lambda_1, \ldots, \lambda_{j_n} )$ as any two such are at least $\beta_0^{j_n}$ apart. Further, being a subset of $F^*$ it is completely contained inside some $B ( \lambda_1, \ldots, \lambda_{j_n} )$. The lat{t}er itself is mapped by $\psi$ into the interval of length $N_{\C} (\beta_0)^{-j_n}$ of $I$ consist{i}ng of numbers whose $N_{\C} (\beta_0)$-adic expansion begins with $0.\lambda_1 \cdots \lambda_{j_n}$. We conclude that $\abs{\psi ( U_n )} \leq N_{\C} (\beta_0)^{-j_n}$ and thereby,
      \begin{align}\label{E:hdim}
        1 &\leq \sum_{n \in \N} \abs{\psi ( U_n )} \leq \sum_{n \in \N} N_{\C} ( \beta_0 )^{-j_n}\\
          &= \sum_{n \in \N} N_{\C} ( \beta_0 )^{-\big\lfloor \frac{\log ( 2 \abs{U_n} )}{\log \beta_0} \big\rfloor} \leq N_{\C} ( \beta_0)\cdot 2^{ \frac{\log N_{\C} (\beta_0)}{\abs{\log \beta_0}}} \sum_{n \in \N} \abs{U_n}^{ \frac{\log N_{\C} (\beta_0)}{\abs{\log \beta_0}}}.\notag
      \end{align}
      \noindent As $\abs{U_n} \leq \delta < \delta_0 < 1$ for all $n \in \N$ and
      \begin{equation}
        \frac{\log N_{\C} (\beta_0)}{\abs{\log \beta_0}} > \sup_{\delta < \delta_0} \frac{\log f (\delta)}{\log \delta} \geq \frac{\log f (\abs{U_n})}{\log \abs{U_n}}
      \end{equation}
      by our assumpt{i}on, we have that
      \begin{equation}\label{E:poslb}
        \sum_{n \in \N} f (\abs{U_n} ) \geq \big( N_{\C} (\beta_0) \big)^{-1} 2^{\log N_{\C} (\beta_0) / \log \beta_0}
      \end{equation}
      for any arbitrary $\delta$-cover $( U_n )$ of $F^*$. Thus, the inf{i}mum of the sums on the lef{t} side of~\eqref{E:poslb} taken over all $\delta$-coverings of $F^*$ is a posit{i}ve number independent of $\delta$. Let{t}ing $\delta \rightarrow 0$ from the right, we conclude that the $f$-dimensional Hausdorf{f} measure of $F^*$ is strictly posit{i}ve. In part{i}cular, this proves our claim.
    \end{proof}
    \begin{corollary}
      Let $\Sigma$ be the full solenoid over $S^1$. Then, the Hausdorf{f} dimension of any \CAW\ subset of $\Sigma$ is inf{i}nite.
    \end{corollary}
    \begin{proof}
      We know $N_{\C} (\beta)$ rises faster than $\beta^{c / ( \beta \log \beta )}$ as $\beta \rightarrow 0$ for some absolute constant $c > 0$, when $\P$ is the set of all rat{i}onal primes. Take $f$ to be the power funct{i}on $r \mapsto r^n$ for some $n \in \N$. The condit{i}on in Proposit{i}on~\ref{P:hdim} is sat{i}sf{i}ed then and we get that the $n$-dimensional Hausdorf{f} measure of any \CAW\ subset of $\Sigma$ is posit{i}ve. F{i}nally, we let $n \rightarrow \infty$.
    \end{proof}
    The same is true for \CAW\ subsets of $\SP$, when \P\ is an inf{i}nite set consist{i}ng of all primes in some arithmet{i}c progression. It will also be interest{i}ng to study the class of exact dimension funct{i}ons for the spaces \XP\ and their \CAW\ subsets~\citep{Fal90}. For f{i}nite \P, we are sat{i}sf{i}ed with a statement about maximality of Hausdorf{f} dimension.
    \begin{proposition}\label{P:findim}
      Let $\abs{\P} < \infty$ and $F \subset \XP$ be a $\beta$-CAW set. Then,
      \[
        \dim F \geq \frac{\log N_{\C} (\beta)}{\abs{\log \beta}}.
      \]
    \end{proposition}
    \begin{proof}
      The proof is the same as that for Proposit{i}on~\ref{P:hdim} till~\eqref{E:hdim} with $\beta$ replacing $\beta_0$ everywhere.
    \end{proof}
    Once again if we let $\beta \rightarrow 0$ for a \CAW\ set, we get
    \begin{corollary}
      Any \CAW\ subset of $\XP$ with $\abs{\P} = l - 1$ has Hausdorf{f} dimension equal to $l$. In part{i}cular, the collect{i}on of points described in Cororllary~\ref{Th:sol} has full dimension.
    \end{corollary}
    \begin{remark*}
      The proofs of Proposit{i}ons~\ref{P:hdim} and \ref{P:findim} are suggest{i}ve that while there can be Schmidt winning subsets in the example metric space of \citeauthor{KW10}~\citep{KW10} ment{i}oned in \S\,\ref{S:games} which are of Hausdorf{f} dimension zero (in fact, countable), any absolute winning subset thereof shall have to be of full dimension.
    \end{remark*}
    
  \section*{Acknowledgments}  
    \noindent The author thanks Anish Ghosh for sugges{t}{i}ng the problem, lots of advice, encouragement and reading various preliminary versions of the manuscript. I am grateful to Jinpeng An, S.~G.~Dani and Sanju Velani for being a pat{i}ent audience and asking many interest{i}ng quest{i}ons.  The calculat{i}ons in Lemma~\ref{L:prod} together with the reference~\cite{MV07} were explained to me by Divyum Sharma.\\[-0.2cm]
    
    \noindent Parts of this work f{i}rst appeared in a slightly dif{f}erent avatar in the author's PhD thesis submit{t}ed to the Tata Inst{i}tute of Fundamental Research, Mumbai in~2017. For this durat{i}on of our invest{i}gat{i}ons, f{i}nancial support from CSIR, Govt.\ of India under SPM-07/858(0199)/2014-EMR-I is duly acknowledged.

\end{document}